\documentclass[preprint,12pt]{elsarticle}
\usepackage{amsmath, amssymb, latexsym, amscd}

\newdefinition{rmk}{Remark}
\newdefinition{definition}{Definition}

\newtheorem{thm}{Theorem}

\newtheorem{cor}{Corollary}
\newtheorem{proposition}{Proposition}

\newproof{pf}{Proof}
\newproof{pf-cor-knots}{Proof of Corollary \ref{knots}}

\DeclareMathOperator{\dem}{dem}

\DeclareMathOperator{\id}{id}

\DeclareMathOperator{\supp}{supp}

\journal{Topology and its Applications}

\begin{document}

\begin{frontmatter}

\title{On a question of B.J.~Baker and M.~Laidacker
concerning disjoint compacta in $\mathbb R^N$
}
\author{Olga Frolkina
\fnref{fn1}}
\ead{olga.frolkina@math.msu.ru}
\address{Chair of General Topology and Geometry\\
Faculty of Mechanics and Mathematics\\
M.V.~Lomonosov Moscow State University\\
and\\
Moscow Center for Fundamental and Applied Mathematics\\
Leninskie Gory 1, GSP-1\\
Moscow 119991, Russia
}

\fntext[fn1]{
The paper was published with the financial support of the Ministry of Education and Science of the Russian Federation as part of the program of the Moscow Center for Fundamental and Applied Mathematics under the agreement No. 075-15-2022-284.}

\begin{abstract}
We 
describe
wild embeddings
of polyhedra into $\mathbb R^N$
which show that
the answer to the question of
B.J.~Baker--M.~Laidacker (1989)
concerning uncountable families of pairwise disjoint
compacta can be twofold. 
The central idea of our construction is the use
of specific wild Cantor sets, namely,
Antoine--Blankinship--Ivanov necklaces and
Krushkal sticky sets.
Our basic tools are Antoine's methods and Shtan'ko demension theory.
\end{abstract}


\end{frontmatter}

\section{Introduction}

The central theme of this work is \emph{disjoint embeddings}:
we would like to
find uncountably many
copies of a given compactum in $\mathbb R^N$ 
simultaneously,
so that they are mutually exclusive.
This kind of problems goes back to 
the result of
R.L.~Moore (1928):
any family of 
pairwise non-intersecting triodes in the plane 
is at most countable.
Continua
that can be placed in $\mathbb R^2$ in an uncountable number
of disjoint homeomorphic copies
satisfy further restrictions on their topological structure,
but many questions still remain open \cite{RSS}, \cite{Hoehn}.

Placing an
uncountable family of disjoint homeomorphic 
copies of a given compactum in the same ambient space
imposes certain restrictions
not only on the compactum itself, but also
 on the behaviour of the embeddings.
As an example:
concentric spheres in $\mathbb R^3$
of arbitrary radii $r>0$
form a family of cardinality continuum.
In contrast to this,
by results of R.H.~Bing, 
it is impossible to place
an uncountable collection
of pairwise disjoint
\emph{wild} closed surfaces in $\mathbb R^3$ (Definition \ref{tame-P}).
A short sketch of Bing's idea can be found in
\cite[Thm.~3.6.1]{BC}; see also \cite{Frolkina} for detailed references.
For $N\geqslant 5$ there is 
a similar impossibility theorem for wild $(N-1)$-spheres in
$\mathbb R^N$
\cite[Thm.~10.5]{Burgess-75}, \cite[p.~383, Thm.~3C.2]{Daverman}, \cite[Thm. 1, 2]{Bryant}.

This paper
concerns with
those collections of pairwise disjoint continua which
are not only homeomorphic,
but ambiently homeomorphic:

\begin{definition} 
Two subsets $X,X'\subset \mathbb R^N$ 
are \emph{ambiently homeomorphic}
(or \emph{equivalently embedded})
if there exists a homeomorphism $h$
of $\mathbb R^N$ onto itself such that
$h(X)=X'$.
In this case, we will write
$h: (\mathbb R^N,X)\cong (\mathbb R^N,X')$.
\end{definition}

As an example, let us recall the result of J.H.~Roberts:
for a ``nicely embedded'' snake-like 
continuum $K\subset \mathbb R^2$
there is an embedding $F:K\times \mathcal C\to\mathbb R^2$
such that each $F(K\times \{t\})$, where $t\in\mathcal C$,
is embedded equivalently to
the given $K$
\cite[Thm.~1]{Roberts}.
(Here $\mathcal C$ is the Cantor set.)
E.K.~van Douwen proved 
that if a separable completely metrizable space~$Y$ 
contains uncountably many mutually disjoint
homeomorphic copies of a given compactum $X$,
then $Y$ contains the product 
$X\times\mathcal C$
\cite[Thm.~1]{vanDouwen};
see 
\cite[Thm.~2.2]{BEM}, \cite[Thm.~1]{Todorcevic}
for different proofs and further improvements.
In general, we can not assert that
the resulting images of $X \times \{t \}$,
for $t\in\mathcal C$,
are embedded equivalently to each other or
to some give copy of $X$.

By
the classical
Lefschetz--Menger--N{\"o}beling--Pontryagin--Tolstowa
Embedding Theorem,
each $k$-dimensional compactum embeds into $\mathbb R^{2k+1}$.
In the case of a $k$-dimensional polyhedron $X$,
the product $X\times I$ embeds into $\mathbb R^{2k+1}$ \cite[Thm. 1.5]{RSS2};
but this statement does not ``see'' how ``individual''
images of $X\times \{t\}$, $t\in I$, are embedded.

\subsection{The question of B.J.~Baker and M.~Laidacker}\label{question}

B.J.~Baker and M.~Laidacker posed
the question \cite[p. 209]{BL}:
\emph{
Let $X\subset \mathbb R^{2k+1}$ be a $k$-dimensional continuum;
is it true that
$\mathbb R^{2k+1}$ 
contains a family 
of pairwise disjoint compacta
$\{ X_\alpha \ | \ \alpha \in \mathcal A \}$, where
$|\mathcal A | = \mathfrak c$
and
$(\mathbb R^{2k+1}, X_\alpha ) \cong ( \mathbb R^{2k+1} , X)$
for each $\alpha \in\mathcal A$~?
}

Corollary~2 of \cite{BL}
gives a positive answer
under an additional ``niceness'' assumption on the embedding of $X$; 
in general, 
the question remained open.
In Section~\ref{statements}, 
\emph{we present two series of examples 
showing that
the answer can be twofold},
depending on more subtle properties of the given
embedding 
$X\subset \mathbb R^{2k+1}$.

In order to clarify which part of the Baker--Laidacker question is open, we need further preparations.

\subsection{On (ambiently) universal spaces}\label{Univ}

In 1916, W.~Sierpi{\'n}ski described a curve
which is now well-known as the Sierpi{\'n}ski carpet, and 
proved that each $1$-dimensional planar compactum 
can be embedded in the carpet
(by that reason, the term ``the Sierpi{\'n}ski 
universal plane curve'' is also used).
In 1921, he proved that the Cantor set 
is a universal space for the class of all zero-dimensional metrizable compacta.

In 1926, K.~Menger
defined a $k$-dimensional generalization $M^k_N\subset \mathbb R^N$
of both the Cantor set and the Sierpi{\'n}ski carpet.
The construction is recalled in Section \ref{M}.
Menger conjectured that
any $k$-dimensional compactum embeds into
$\mathbb R^{2k+1}$, and proved it for $k=1$.
For arbitrary $k$ this was proven
in 1931 independently
by G.~N{\"o}beling;
by L.S.~Pontrjagin and G.~Tolstowa;
and
by S.~Lefschetz.
Menger conjectured
that $M^k_N$ is a universal space for the class of all
$k$-dimensional compact subsets of $\mathbb R^N$, and
proved it for two cases: $k=1$, $N=3$; and $k=N-1$ (for references, see \cite{Chigogidze}; the pictures of
the Menger universal curve $M_3^1$ 
also known as the Menger sponge, the Menger cube, 
the Sierpi{\'n}ski cube, or the Sierpi{\'n}ski sponge, can be found in \cite[p.~131--132]{DV}).

Embeddability of any
$k$-dimensional compact metric space 
into $M^k_{2k+1}$ was shown by S.~Lefschetz \cite{Lefschetz}.
In its full generality, 
Menger's conjecture
was proven by
M.A.~Shtan'ko:
each $k$-dimensional compactum $X$ embeddable in $\mathbb R^N$
can be embedded into $M^k_N$
\cite[Thm.~1]{Stanko1971}; see also
\cite{Shtanko-diss}
or
\cite[Cor.~5.5.4]{DV}.
But if $X$ is already embedded in $\mathbb R^N$, 
we can not guarantee
that 
$X$ can be position-wise embedded in $M^k_N$
in the following sense:

\begin{definition}
Let $X,Z\subset \mathbb R^N$.
We say that \emph{$X$ can be ambiently embedded
in $Z$},
or 
\emph{$X$ can be position-wise embedded in~$Z$} \cite{BL},
if there exists a homeomorphism $h$ of $\mathbb R^N$ onto itself such that $h(X)\subset Z$.
\end{definition}

For $(k,N)=(0,1)$ and $(k,N)=(1,2)$
the existence of an ambient homeomorphism 
$h: \mathbb R^N\cong \mathbb R^N$
that takes a given $k$-dimensional compactum $X$ 
to $h(X) \subset M^k_N$
follows from arguments of Sierpi{\'n}sky,
and for $(k,N)=(0,2)$ --- 
from Antoine's result
\cite[Cor. 2.3.2]{Keldysh}, \cite[Thm. 13.7]{Moise}.
In general,
the existence of such an $h$
is equivalent to the inequality
$\dem X \leqslant k$
\cite[Thm.~2]{Stanko1971},
\cite{Shtanko-diss}, 
\cite[Thm.~3.5.1]{DV}
(partial results were also obtained in
\cite{Bothe1964univ}, \cite{Bryant1969}).
Here ``$\dem $'' denotes 
``demension''.
This word does not contain mistakes:
it abbreviates ``dimension of embedding''
introduced by M.A.~Shtan'ko.
The inequality $\dem X \leqslant k$
means that $X\subset \mathbb R^N$ behaves geometrically
much like a polyhedron of dimension $\leqslant k$.
As an example, $\dem M^k_N = k$.

Shtan'ko demension theory
was developed with the aim of extending 
the notion of tameness from polyhedra 
(Section \ref{wildness})
to arbitrary compacta,
see
\cite{Stanko1970},
\cite{Shtanko-diss}, 
\cite{Edwards}, \cite[3.4, 3.5]{DV},
\cite{Cernavsky}, \cite{Melikhov}.
We will not set forth the foundations of this theory in the present paper, 
indicating the suitable source if necessary.

\subsection{The construction of Menger compacta $M^k_N$}\label{M}

Let $0\leqslant k\leqslant N$.
Take the standard cube $I^N = [0,1]^N\subset \mathbb R^N$.
Let
$\mathcal T_j$
be
the collection of all $N$-cubes
obtained by subdividing $I^N$ into $3^{jN}$
congruent $N$-cubes
by hyperplanes drawn perpendicular to the edges of $I^N$
at points dividing the edges into $3^j$ equal segments.
The compactum $M^k_N$ is constructed inductively.
The family $\mathcal F_0 $ consists of the only element $ I^N$.
The collection
$\mathcal F_{j+1}$
consists of all $v\in \mathcal T_{j+1}$
with the property:
$v$ is a subset of some element $w\in \mathcal T_j$
such that $v$ intersects some $k$-face of~$w$.
Now let $P_j $ 
be the union of all elements of $\mathcal F_j$.
\emph{The $k$-dimensional Menger compactum constucted in $\mathbb R^N$}
is defined by
$$
M^k_N = \bigcap\limits_{j=0}^\infty  P_j .
$$

Taking $k=0$, we get $M^0_N = \mathcal C ^N\cong \mathcal C$,
where $\cong $ stands for a homeomorphism.

For $1\leqslant k\leqslant N$, the compactum $M^k_N$ 
is connected, therefore it is also called
\emph{the $k$-dimensional Menger continuum constucted in $\mathbb R^N$}.
We also say
\emph{the standard $k$-dimensional Menger continuum in
$\mathbb R^N$},
especially when we consider
other
inequivalent embeddings 
$M^k_N\hookrightarrow \mathbb R^N$.
Observe
that
$$
\mathcal C^N = M^0_N \subset M^1_N\subset \ldots \subset M^{N-1}_N \subset
M^N_N = I^N .
$$

For $N\geqslant 2k+1$ we have
$M^k_N\cong M^k_{2k+1}$ \cite{Bestvina}.

Let us emphasize that $M^k_N\subset \mathbb R^N$ is thought as
a specific subset rather than an abstract compactum;
it should be distinguished from
an arbitrary
embedding $M^k_N\hookrightarrow \mathbb R^N$ 
(compare with Proposition \ref{3-tangled}).

\subsection{The case considered by B.J.~Baker and M.~Laidacker}

B.J.~Baker and M.~Laidacker 
constructed, 
for any $k\geqslant 0$,
an embedding
$\Xi : M_{2k+1}^k \times\mathcal C 
\rightarrow \mathbb R^{2k+1}$
such that 
all compacta
$\Xi ( M_{2k+1}^k \times \{ t\})$, where $t\in\mathcal C$,
are obtained from each other by parallel translations.
Moreover, by
\cite[Prop. 2]{Stanko1971}
(see also \cite[Prop. 9]{Shtanko-diss} or \cite[Prop. 3.5.2]{DV})
each of $\Xi ( M_{2k+1}^k \times \{ t\})$ is embedded equivalently
to the standard Menger compactum 
$M_{2k+1}^k \subset \mathbb R^{2k+1}$
described in Section~\ref{M}.
Hence their question (Section \ref{question}) has a positive answer
for those $X\subset \mathbb R^{2k+1}$ that satisfy $\dem X \leqslant k$
(equivalently, for $X$ ambiently embeddable into $M^k_{2k+1}$)
\cite[Cor.~2]{BL}.
Our paper is devoted to those
$X$ that can not be ambiently embedded into $M^k_{2k+1}$.

\begin{rmk}
\cite[Cor. 2]{BL}
strengthens
the classical
Lefschetz--Menger--N{\"o}beling--Pontryagin--Tolstowa
theorem
since each
$k$-dimensional compactum can be embedded into
the standard Menger compactum
$M^k_{2k+1}\subset \mathbb R^{2k+1}$
\cite{Lefschetz}.
\end{rmk}

\subsection{On wild embeddings}\label{wildness}

Embeddings in our examples are wild.
Theory of wild embeddings 
appeared at the beginning of the 20th century,
in attempts to rigorously prove the
Schoenflies theorem 
\cite[Thm. 9.6, 10.3]{Moise}
and to generalize it to higher dimensions.
Let us recall

\begin{definition}\label{tame-C}
A zero-dimensional compact 
set $X\subset \mathbb R^N$ is called \emph{tame}
if there exists
a homeomorphism $h$ of $\mathbb R^N$ 
onto itself such that 
$h(X)$ lies on a straight line;
otherwise, $X$ is called \emph{wild}.
\end{definition}

\begin{definition}\label{tame-P}
A subset of $\mathbb R^N$
is called a \emph{polyhedron}
if it is the union of a finite collection of simplices.
A compactum $X \subset\mathbb  R^N$  
homeomorphic to a polyhedron is
called
\emph{tame}
if there exists a homeomorphism  
$h$ of $\mathbb R^N$ onto itself 
such that $h(X)$ is a polyhedron in $\mathbb R^N$;
otherwise, $X$ is \emph{wild}.
\end{definition}

In 1921, L.~Antoine proved that
 each zero-dimensional compactum
 in $\mathbb R^2$
is tame.
In 1920 he sketched and in 1921 explicitly constructed
a first wild Cantor set in $\mathbb R^3$
now called Antoine's necklace \cite[Section 18]{Moise}.
His construction was extended to higher dimensions
independently by A.A.~Ivanov 
\cite{Ivanov},
\cite{Ivanov-diss}
and by
W.A.~Blankinship
\cite{Blankinship};
see also
\cite{Eaton} and \cite[Example 4.7.1]{DV}.
Besides
widely known Antoine's necklace,
there is an essentially different construction
of a wild Cantor set in $\mathbb R^3$ 
given 
by P.S.~Urysohn in 1922--23; for references, see \cite{Frolkina}.
Now, examples of wild Cantor sets in $\mathbb R^N$ include the
ones with such strong properties as
simply-connectedness of the complement,
slipperiness,
stickiness (see Section \ref{On-Cantor-sets}).
The behaviour of wild Cantor sets resembles that
of polyhedra of codimension $2$
\cite[Thm. 1.4]{Edwards}, \cite[Thm. 3.4.11]{DV}.

Examples of wild $2$-spheres in $\mathbb R^3$ (L.~Antoine 1921, J.W.~Alexander 1924), 
wild arcs in $\mathbb R^3$ (L.~Antoine 1921; R.H.~Fox and E.~Artin 1948)
and even
everywhere wild arcs in $\mathbb R^3$ (L.~Antoine 1924)
are 
well-known.
For $N\geqslant 3$,
each uncountable compact
subset of $S^{N-1}$
can be embedded in $\mathbb R^N$ in uncountably many
inequivalent wild ways (R.B.~Sher 1968, D.G.~Wright 1986).
For further information on topological embeddings, refer to
\cite{Keldysh},
\cite{Rushing},
\cite{Moise}, \cite{DV},
\cite{Chernavskii-review},
\cite{ZC-review},
\cite{BC},
\cite{Burgess-75}, \cite{Daverman}.
For embeddings of arbitrary compact sets, the concept of tameness was introduced and deeply studied by M.A.~Shtan'ko 
\cite{Stanko1970}--\cite{Shtanko-diss}. 
See Section \ref{Univ} and references therein.

Our examples 
are based on Antoine's ideas and methods,
and on results of Shtan'ko.

\subsection{Notation and agreements}

All spaces are supposed to be metrizable, all maps are continuous.
As a rule, the metric is denoted by $d$, 
this does not lead to ambiguous understanding.

$\mathbb R^N$ denotes Euclidean $N$-dimensional space
with the usual metric.
$S^N$ is the standard unit sphere in $\mathbb R^{N+1}$ with the induced metric.
$I=[0,1]$.

The symbol $\mathfrak X$
is used for an abstract space,
to distinguish it from its 
homeomorphic copies which are embedded into
a space $Y$;
such embedded copies of $\mathfrak X$
are denoted by $X$ or $X_\alpha $.

For a subset $A$ of a metric space $Y$,
the symbol $\overline A$ denotes the closure,
$\mathring A$ the interior,
and $O_{\varepsilon }(A) = \{ y\in Y \ |\  d(y,A) < \varepsilon \}$
the $\varepsilon $-neighbourhood.
For non-empty compact subsets $A,B\subset Y$,
the distance $d(A,B)$ is the minimum of distances
between $a\in A$ and $b\in B$.

The symbol $Y\cong Y$ denotes a homeomorphism of $Y$ onto itself.
Similarly, 
$h: (Y , A)\cong (Y , A)$ is a homeomorphism 
of pairs, i.e. a homeomorphism
$h: Y\cong Y$ such that
$h(A)= A$.

A self-homeomorphism $h$ of a metric space $Y$ is called
an $\varepsilon $-homeomorphism
if $d(x, h(x))\leqslant \varepsilon  $
for each $x\in Y$.

For a homeomorphism $h: Y\cong Y$,
its support $\supp h$
is the closure of the set
$\{ x\in Y \ | \ h(x) \neq x \}$.

As usual, $\id $ is the identity map.

For a topological manifold $M$,
 its boundary is  denoted by
$\partial  M$.

$\mathcal C$ is the Cantor set.

$| \cdot |$ is the cardinality of a set;
$\mathfrak c = | \mathbb R | $
is the cardinality of the continuum.

$\mathcal A$ is usually an index set.

\section{Statements}\label{statements}

\subsection{Negative examples}

In our first ``negative'' result,
the dimension of the ambient
space is at least~$4$.
The $3$-dimensional case is treated below, in Proposition~\ref{3-tangled}.

\begin{thm}\label{BL-negative-m}
For any $N\geqslant 4$
and any uncountable compactum $\mathfrak X\subset S^{N-1}$
there exists an embedding
 $f :\mathfrak X\rightarrow \mathbb R^N$ such that
 \\
 1)
the image $X:=f( \mathfrak X)$ can not be position-wise embedded in $M^{N-3}_N$,
and
\\
2) it is impossible to place in $\mathbb R^N$
uncountably many
pairwise disjoint compacta ambiently homeomorphic to $X$.
\end{thm}

Theorem \ref{BL-negative-m} is proved in Section \ref{proof-BL-negative-m}.

As a consequence, we get ``negative'' examples
for the Baker--Laidacker problem, assuming that the dimension
of the ambient space is at least $5$:

\begin{cor}\label{BL-negative}
For each $k\geqslant 2$
there exists an embedding 
$f :I^k\rightarrow \mathbb R^{2k+1}$
such that
\\
1) the $k$-cell $Q:=f (I^k)$ 
can not be position-wise embedded into $M^k_{2k+1}$
(and even into $M^{2k-2}_{2k+1}$), and
\\
2) it is impossible to place in $\mathbb R^{2k+1}$ 
uncountably many pairwise disjoint $k$-cells
embedded equivalently to $Q$.
\end{cor}

Corollary \ref{BL-negative}
can not be ``word-to-word'' extended to the case of 
$2k+1=3$,
since each
arc in $\mathbb R^3$
can be position-wise embedded
in $M^1_3$
\cite[Satz~4]{Bothe1964univ}.
That is, desired examples should have a more complicated
structure.

\begin{proposition}\label{3-tangled}
There exist an
embedding
$f: M^1_3 \rightarrow \mathbb R^3$
and a closed broken line $L\subset \mathbb R^3$
such that the compactum 
$X=L\cup f(M^1_3)$
has the following properties:
\\
1) $X$ is connected,
\\
2) $\dim X = 1$,
\\
3) $X$ can not be position-wise embedded into $M_3^1$,
\\
4)  it is impossible to place in $\mathbb R^3$ 
uncountably many pairwise disjoint compacta
embedded equivalently to $X$.
\end{proposition}

\begin{pf}
The construction of
the standard  Menger continuum
$M^1_3\subset \mathbb R^3$ 
can be thought of as starting from $I^3$
and drilling through rectangular channels.
If we drill knotted channels instead,
taking care of their size and position,
the resulting compactum $T\subset \mathbb R^3$
will be homeomorphic
to $M^1_3$ by the Anderson Characterization Theorem 
\cite[Thm. XII]{Anderson};
this is the image of the desired embedding: $T=f(M^1_3)$.
And it will be impossible
to remove some
unknotted closed broken line $L$
from $T$ by a small push.
With details and proof,
such an embedding is presented in
\cite{Bothe1964};
with less details, 
a similar example is described 
 in
\cite{McMillan-Row}.
Without proof, a similar
construction was given
in
\cite[p. 788--789]{Frankl-Pontrjagin} 
in connection with the problem of
coincidence of 
Brouwer--Menger--Urysohn and Alexandroff
dimensions.
Thereby,
we will refer to the paper
\cite{Bothe1964} containing maximum details.
The desired properties follow  from Bothe's results.
Indeed, 1) is by construction. 2) is proved in \cite{Bothe1964}.
The set
$T$ can not be position-wise embedded
in $M^1_3$ by \cite{Bothe1964}, this implies 3).
Finally, let us prove 4). Suppose the contrary; by Corollary
\ref{R-uncountably},
for any $\varepsilon >0$
there is a homeomorphism 
$g:\mathbb R^3\cong \mathbb R^3$ such that
$ X \cap  g(X) = \emptyset $ and
$d(x,g(x)) \leqslant \varepsilon $ for each $x\in X$.
Hence $T\cap g(L) = \emptyset $;
this is contradicted by \cite[p. 255]{Bothe1964}.
\end{pf}

\begin{rmk}
In Proposition \ref{3-tangled}, 
the compactum $T=f(M^1_3)$ itself already 
has properties 1)--3). 
We conjecture that $T$ also has property 4);
verification of this would require additional
reasoning which, in our opinion, can be carried out similarly
to \cite{Bothe1964}.
\end{rmk}

\subsection{Affirmative examples}

\begin{thm}\label{BL-bouquet}
Suppose that $N\geqslant 4$,
$1\leqslant k\leqslant N-3$, and
$2k+1 \leqslant N$.
Let $P\subset \mathbb R^N$ be a $k$-dimensional
polyhedron
such that
for some point $a\in P$,
its neighbourhood in $P$
is homeomorphic to the interval $(0,1)$.
Then 
there exists an embedding
$F: P \times \mathcal C \rightarrow \mathbb R^N$
with the properties:
\\
1) all compacta $F( P \times \{ t\})$, where $t\in\mathcal C$,
are embedded into $\mathbb R^N$ equivalently to each other, and
\\
2) 
no one of $F( P \times \{ t\})$, where $t\in\mathcal C$,
can be position-wise embedded into
$M^{N-3}_N$ (hence also into $M^k_N$).
\end{thm}

For the proof, see Section \ref{proof-BL-bouquet}.
As a consequence,
we get a series of ``affirmative'' examples 
for the Baker--Laidacker question:

\begin{cor}\label{BL-bouquet-cor}
Let $k\geqslant 2$.
Let $P$ be a $k$-dimensional compact
polyhedron
$P\subset \mathbb R^{2k+1}$ 
such that
for some point $a\in P$,
its neighbourhood in $P$
is homeomorhic to the interval $(0,1)$.
Then
there exists an embedding
$F: P\times \mathcal C \rightarrow \mathbb R^{2k+1}$
such that the compacta 
$F( P\times \{ t\})$, $t\in\mathcal C$,
are pairwise 
ambiently homeomorphic,
and no one of them can be position-wise embedded into 
$M^k_{2k+1}$ (and even into $M^{2k-2}_{2k+1}$).
\end{cor}

\begin{rmk}
In \cite{Frolkina},
 for each 
$N\geqslant 3$ and each 
$1\leqslant k \leqslant N-1$
we described
an embedding 
$F: I^k\times \mathcal C \rightarrow \mathbb R^{N}$
such that all $k$-cells
$F( I^k\times \{ t\})$, $t\in\mathcal C$,
are \emph{inequivalently} embedded in
$\mathbb R^N$; no one of these cells can be position-wise embedded 
in $M^{N-3}_{N}$.
(For $N=2$, the situation is different:
any Cantor fence $I\times \mathcal C
\hookrightarrow \mathbb R^2$
is ambiently equivalent to the standard one \cite[Thm. 2]{B1958}, \cite{TW}.)
We believe that both Theorem~\ref{BL-bouquet}
and Corollary~\ref{BL-bouquet-cor} hold true
for $P=I^k$;
this probably can be obtained
using ideas from \cite{Frolkina}.
\end{rmk}

Recall that the question of Baker and Laidacker
concerns connected spaces (continua).
But it also makes sense for non-connected spaces.
In this case, another construction is possible:

\begin{thm}\label{BL-disconnected}
For each $N\geqslant 3$
and each $0\leqslant k\leqslant N-3$
there exists an embedding
$F: (I^k\times \mathcal C)\times \mathcal C \rightarrow \mathbb R^N$
such that
\\
1) all compacta
$F( (I^k\times \mathcal C)\times \{ t\})$, for $t\in\mathcal C$,
are embedded in $\mathbb R^{N}$
equivalently to each other,
and 
\\
2) no one of these compacta can be position-wise embedded into
$M^{N-3}_{N}$ (hence into $M^k_N$ also).
\end{thm}

Taking $k\geqslant 2$ 
and $N=2k+1$
in Theorem \ref{BL-disconnected},
we get new affirmative
examples for the Baker--Laidacker problem,
except for the connectivity requirement.
The proof of Theorem \ref{BL-disconnected}
is given in Section \ref{proof-BL-disconnected}.

\section{Proofs}

\subsection{Main observations regarding disjoint embeddings}

In our proofs we rely on the following assertion.
It is probably well known (at least partially), but
we did not find it in the literature.
The proof does not contain 
fundamentally 
new ideas;
however, we present it 
for completeness.
Condition (iii) is rather strong,
as is shown in \cite{CW} and 
in a recent arxiv preprint
``An answer to a question of J.W.~Cannon and S.G.~Wayment''
 by the author.

\begin{thm}\label{main}
Let $\mathfrak X$, $Y$ be compact spaces.
Let $e : \mathfrak X \to Y$ be an embedding;
denote $X=e (\mathfrak X)$.
The following are equivalent:
\\
(i) there exists an embeding $\Psi : \mathfrak X \times \mathcal C \to Y$
such that
$(Y, \Psi (\mathfrak X \times \{ t \} )) \cong (Y,X)$ 
for each $t\in\mathcal C$,
\\
(ii) 
there exists a collection 
$\{ X_\alpha \ |  \ \alpha \in \mathcal A\}$
of 
$|\mathcal A| = \mathfrak c$
pairwise disjoint subsets of $Y$
such that 
$(Y, X_\alpha )\cong (Y,X)$ for any $\alpha \in\mathcal A$,
\\
(iii) for any $\varepsilon >0$ there exists a homeomorphism
$h: Y\cong Y$ such that
$X\cap h( X) = \emptyset $ and
$d(x, h(x))\leqslant \varepsilon $ for each $x\in Y$.
\end{thm}

\begin{pf}
Let $d$ be any fixed metric on $Y$,
and endow
the space
$C(Y,Y)$
of continuous maps
with the distance
$\rho (f,g) = \sup \{ d (f(x) , g(x)) \ | \ x\in Y\}$.

$(iii)\Rightarrow (i)$
is proved by the standard Cantor tree type argument,
similarly to \cite{Roberts},
\cite{Bing1962};
compare with
 \cite[Exercise 4.8.1 and p.~180, Remark]{DV},
\cite{Eaton}.

By compactness of $Y$,
property $(iii)$
remains valid after replacing $X$ by any $X'$
such that
$(Y,X')\cong (Y,X)$.

Let $\{\varepsilon _m\}$ be a sequence of positive numbers.
(In fact, each of its elements will be selected based on what has already been done in the previous step; the details will become clear later.)

\emph{Step 1.}
Apply $(iii)$ to $X$.
Take homeomorphisms $h_0, h_1 : Y\cong Y$ such that
$\rho (h_0, \id _Y) \leqslant \varepsilon _1 $,
$\rho (h_1, \id _Y) \leqslant \varepsilon _1 $,
$h_0 ( X) \cap h_1 (X)  = \emptyset $.

\emph{Step 2.}
Apply $(iii)$ to $h_0(X)$ and to $h_1(X)$.
Take homeomorphisms 
$h_{00}, h_{01}, h_{10}, h_{11}  : Y\cong Y$
such that
$$
\max \{ \rho (h_{00}, \id_Y),
\
\rho (h_{01}, \id_Y) ,
\
\rho (h_{10}, \id_Y) ,
\
\rho (h_{11}, \id_Y) 
\}
 \leqslant \varepsilon _2,
$$
and 
$$h_{00} h_0 ( X) ,\ 
 h_{01} h_0 (X)  ,\
 h_{10} h_1 (X) ,\
 h_{11} h_1 (X)  $$
 are pairwise disjoint.

We continue in the same way.
On \emph{Step $m$} we 
apply $(iii)$ to $2^{m-1}$
compacta constructed on \emph{Step $(m-1)$}, and
obtain
a collection of homeomorphisms
$\{ h_{a_1a_2\ldots a_m} \ | \ a_1,\ldots , a_m \in \{ 0;1\}\}$
such that 
$$
\max 
\{  \rho (h_{a_1a_2\ldots a_m} , \id _Y) \ | \  a_1,\ldots , a_m \in \{ 0;1\} \} \leqslant \varepsilon _m
$$
and
all $2^m$ compacta 
$$
h_{a_1a_2\ldots  a_m} h_{a_1a_2\ldots a_{m-1}} \ldots h_{a_1a_2} h_{a_1} (X) 
\quad\text{ for }
\quad
a_1,\ldots , a_m \in \{ 0;1\}
$$
are pairwise disjoint.
Assume that
the sequence $\{ \varepsilon _m\}$
is
rapidly decreasing,
so that we have:
for each 
$\sigma = (a_1,a_2,a_3, \ldots ) \in \{0,1\} ^\mathbb N$
the sequence
$$h_{a_1} ,\  h_{a_1 a_2} h_{a_1}, \ h_{a_1 a_2 a_3} h_{a_1 a_2} h_{a_1} , 
\ldots $$
converges to a homeomorphism
$h_{\sigma } : Y\cong Y$
\cite{Roberts},
\cite[Thm.~7]{Bing1962},
\cite[Prop.~2.2.2]{DV}.

Fix any homeomorphism
$\xi : \mathcal C \cong \{0;1\}^\mathbb N$.
The desired embedding $\Psi $ is defined by
$$
\mathfrak X\times \mathcal C 
\stackrel{e\times \xi }{\cong } X \times  \{0;1\}^\mathbb N
\rightarrow Y,
\quad
(\mathfrak x, t)
\mapsto 
(e (\mathfrak x),  \xi (t) )
\mapsto
h_{\xi (t)}  (e (\mathfrak x)) .
$$

For $t,t'\in \mathcal C$ 
the desired homeomorphism
$(Y, \Psi ( \mathfrak X \times \{t\} )) \cong 
(Y , \Psi( \mathfrak X \times \{t' \} ))$
can be defined as the composition 
$h_{\xi (t')} (h_{\xi (t)})^{-1} $.

\medskip

$(i)\Rightarrow (ii)$ is evident.

\medskip

$(ii)\Rightarrow (iii)$.
For each $\alpha \in \mathcal A$ fix a homeomorphism
$g_\alpha : (Y , X) \cong (Y, X_\alpha )$.
The space 
of continuous maps
$C (Y , Y)$ 
is complete by compactness of $Y$.
Hence
the uncountable set
$ \{ g_\alpha \ | \  \alpha \in \mathcal A\}$
contains a countable subset of pairwise distinct elements
$f_{\star }, f_1, f_2,\ldots  $ such that
$\lim\limits _{m\to\infty } \rho (f_m , f_\star ) =0$.
Let $h_m :=  (f_\star )^{-1} f_m  : Y\cong Y$.
We have
$$
X\cap h_m (X)  = (f_\star )^{-1} (f_\star (X) \cap f_m (X))  =
(f_\star )^{-1} (\emptyset ) = \emptyset .
$$
It is easy to verify that $
\lim\limits _{m\to\infty }  \rho ( \id _{Y} , h_m  )=0$,
hence $X$ satisfies $(iii)$.
\end{pf}

\begin{rmk}\label{orient}
Analyzing the proof of $(ii)\Rightarrow (iii)$ of Theorem \ref{main},
we easily see that in the case of $Y=S^N$ 
(or any oriented closed manifold),
$h$ can be taken to be orientation-preserving.
\end{rmk}

\begin{rmk}\label{vD}
In the case of $Y=\mathbb R^N$, conditions 
$(i)$ and $(ii)$
of Theorem \ref{main}
are equivalent.
This can be shown following the same proof,
applied to the one-point compactification $S^N$;
all maps $g_\alpha $ and $h_{a_1\ldots a_m}$
should be taken to preserve
the added compactifying point.
\end{rmk}

\begin{cor}\label{push}
Let $\mathfrak X$ be a compact space.
Let $e : \mathfrak X \to \mathbb R^N$ be an embedding, and
denote $X=e (\mathfrak X)$.
Suppose that for any 
$\varepsilon >0$ 
there exists a homeomorphism
$h: \mathbb R^N\cong \mathbb R^N$
with the properties:
\\
1) $X\cap h(X) = \emptyset $,
\\
2) $\supp h \subset O_{\varepsilon }(X) $,
and
\\
3) $d(x,h(x))\leqslant \varepsilon $ for each
$x\in \mathbb R^N$.

Then there exists an embedding
$\psi : \mathfrak X \times \mathcal C \to \mathbb R^N$
such that
$(\mathbb R^N, \psi (\mathfrak X \times \{ t \} )) \cong (\mathbb R^N, X)$
for each $t\in\mathcal C$.
\end{cor}

\begin{pf}
Let $Y$ be
a closed ball in
$\mathbb R^N$ that contains $X$ in its interior.
By assumptions,
for each
$\varepsilon >0$ 
there exists an $\varepsilon $-homeomorphism
$f: Y\cong  Y$ such that
$X\cap f(X) = \emptyset $
and $f|_{\partial Y} = \id $.
By
Theorem~\ref{main} 
we get an embedding
$g: \mathfrak X \times \mathcal C \to  Y$ 
such that
$(Y, g(\mathfrak X \times \{ t \} )) \cong (Y,X)$
for each $t\in\mathcal C$.
The desired embedding $\psi $
can be obtained as the composition of $g$ and the
inclusion $Y \subset \mathbb R^N$.
\end{pf}

\begin{cor}\label{R-uncountably}
Let $X\subset \mathbb R^N$ be a compact set.
Suppose that there exists
an uncountable collection
$\{ X_\alpha \ | \ \alpha \in\mathcal A\}$
of subsets of $\mathbb R^N$ 
such that
$(\mathbb R^N, X_\alpha )\cong (\mathbb R^N,X)$ for each 
$\alpha \in \mathcal A$.

Then for each bounded open subset $U\subset \mathbb R^N$
containing $X$
and for each 
$\varepsilon >0$
there exists an orientation-preserving homeomorphism 
$g:\mathbb R^N\cong \mathbb R^N$
such that
\\
1) 
$X\cap g(X) = \emptyset $,
and
\\
2) $d(x, g(x)) \leqslant \varepsilon $ for any $x\in U$.
\end{cor}

Observe that we do not require $g(U)$ to be a subset of $U$.

\begin{pf}
Let $ q\in S^N$ be the north pole. 
Let
$p : S^N\setminus \{ q\} \to \mathbb R^N$ be the stereographic projection map,
and
$i= p^{-1}: \mathbb R^N \rightarrow S^N$.
In this proof, 
$d$ denotes the distance in $\mathbb R^N$; 
the standard distance in $S^N$ 
induced from $\mathbb R^{N+1}$ is denoted by 
$D$.

Suppose that $U$ and $\varepsilon $ are given.
Take a bounded open subset
$V\subset \mathbb R^N$ such that
$\overline U \subset V$.

For brevity, let
$\widehat X := i(X)$, 
$\widehat X_\alpha := i(X_\alpha )$,
$\widehat U := i(U)$,
$\widehat V := i(V)$.
Observe that $(S^N , \widehat X_\alpha )\cong (S^N, \widehat X)$ for each $\alpha \in\mathcal A$.

The closure of $\widehat V $ in $S^N$
is a compact set which does not contain~$q$.
Hence the restriction map
$p|_{\widehat V} : \widehat V \to V$
is well-defined and uniformly continuous. 
Choose a $\gamma >0$ such that
for any $ x,  x' \in \widehat V$ with
$D( x ,  x')\leqslant \gamma $,
we have
$d(p  (x), p (x'))\leqslant \varepsilon $.
We may moreover assume that
$$
O_{\gamma } (\widehat U ) \subset \widehat V
\quad
\text{and}
\quad
O_{\gamma } (\widehat U ) \cap O_{\gamma } (q) = \emptyset .
$$
Using Remark \ref{orient},
take an orientation-preserving homeomorphism
$h:S^N\cong S^N$ such that
${\widehat X \cap h(\widehat X) = \emptyset }$
and
$D ( x , h ( x)) \leqslant \frac{\gamma }{2}$
for each $ x \in S^N$.
In particular, we have $D ( q , h ( q)) \leqslant \frac{\gamma }{2}$. Take
an orientation-preserving homeomorphism
$f : S^N\cong S^N$ such that
$f  h(q) = q$ 
and
$\supp f \subset O_{ \gamma } (q)$.
Let 
$h' : 
S^N\setminus \{q\} \cong 
S^N\setminus \{q\}  $ denote the restriction
of $fh : S^N\cong S^N$.

Define
the desired self-homeomorphism of $\mathbb R^N$
by the formula 
$g=p  h'  i$.
To verify 1), observe that
$\widehat U \cup h(\widehat U) \subset S^N\setminus O_\gamma (q) $,
hence 
$h'|_{\widehat U}= h|_{\widehat U}$, and we get
$$
X\cap g(X) = 
p (\widehat X )\cap ph' (\widehat X ) =
p( \widehat X \cap h'(\widehat X ))=
p( \widehat X \cap h(\widehat X ))
= \emptyset .
$$
To prove 2), 
take any 
$x\in U$ and
denote 
$\widehat x := i(x)\in\widehat U$.
Recall that
$$
D(\widehat x, h' (\widehat x)) = 
D( \widehat x, h(\widehat x)) \leqslant\frac{\gamma }{2} .
$$
Consequently $\widehat x$ and 
$h(\widehat x )$
lie in
$O_\gamma ( \widehat U)\subset  \widehat V$.
Together with the choice of $\gamma $, this implies
$$
d(x,g(x)) = d(p(\widehat x), ph'(\widehat x)) \leqslant \varepsilon .
$$
Corollary \ref {R-uncountably} is proved.
\end{pf}

\subsection{Pushing Cantor sets off themselves}\label{On-Cantor-sets}

R.J.~Daverman conjectured
\cite[Conj. 1]{Daverman-question}
that
for any two Cantor sets $X$ and $X'$ in $\mathbb R^N$
and any
$\varepsilon >0$
there is an $\varepsilon $-homeomorphism $h:\mathbb R^N \cong \mathbb R^N$
such that
$X \cap h(X') = \emptyset $.
This is known to be true for $N\leqslant 3$
\cite[Thm.~1]{Sher1969}.

However, for each $N\geqslant 4$
V.~Krushkal constructed
a sticky (wild) Cantor set in $\mathbb R^N$
\cite[Thm. 1.1]{Krushkal}:
it cannot be
isotoped off of itself by any sufficiently small 
ambient isotopy. 
(Compare \cite{Wright-pushing}.)

Modifying Alexander's idea,
J.~Kister observed 
that
each pair of 
 $\varepsilon $-close
 homeomorphisms of
$\mathbb R^N$
is connected by an ambient $\varepsilon $-isotopy
\cite[Thm.~1]{Kister}
(for general manifolds, see \cite{Chernavskii1969}, \cite{Chernavskii2008}).
Hence it is impossible to take the Krushkal set
off itself by any sufficiently small 
ambient homeomorphism.
This can also be derived directly from Krushkal's
arguments.
In fact, Krushkal
constructed a Cantor set
$K\subset \mathbb R^N$
together with a bounded open neighbourhood $U\subset \mathbb R^N$
and
an
$\varepsilon >0$
so that the following is satisfied:
if $h:\mathbb R^N\cong\mathbb  R^N$
is an orientation-preserving homeomorphism
satisfying
$d(x,h(x))\leqslant \varepsilon $ 
for each
$x\in U$, then $K\cap h(K)\neq\emptyset $.
 Corollary~\ref{R-uncountably} now implies:

\begin{cor}\label{Krushkal-Cantor}
Let $K\subset\mathbb R^N$ be a Krushkal Cantor set,
$N\geqslant 4$.
Suppose that
$\{ X_\alpha  \ | \ \alpha \in \mathcal A \}$ is a family 
of
pairwise disjoint Cantor sets in $\mathbb R^N$ such that
each $X_\alpha $ is
ambiently homeomorphic to $K$.
Then
$\mathcal A$ is no more than countable.
\end{cor}

This Corollary will be used in the proof of Theorem~\ref{BL-negative-m}.

\begin{rmk}
The property mentioned immediately before Corollary \ref{Krushkal-Cantor}
implies: 
if $K\subset \mathbb R^N$
is the Cantor set constructed by Krushkal
and
$g:\mathbb R^N\cong\mathbb R^N$
is a homeomorphism,
then $g(K)$ is also a sticky Cantor set.
For this reason,
by the Krushkal Cantor set we will
mean any Cantor set embedded equivalently
to the concrete one described by Krushkal.
\end{rmk}

In contrast to Kruskal Cantor sets,
it is possible to place in $\mathbb R^N$
a family of $\mathfrak c$
pairwise disjoint equivalently
 embedded  Antoine--Blankinship--Ivanov Cantor sets
(see Section \ref{wildness}).
This will be used in the proof of Theorem \ref{BL-bouquet}.

The demension of each wild Cantor set
in $\mathbb R^N$
equals $N-2$
(see
\cite[Thm. 1.4]{Edwards}, \cite[Thm. 3.4.11]{DV};
in the case of $N=4$, refer to
\cite{Shtanko-4dim} or \cite[p.~5]{BDVW}).
The polyhedra constructed in our paper
can not be position-wise embedded 
into $M^{N-3}_N$
since they
contain
wild Cantor sets.
                                
\subsection{Proof of Theorem \ref{BL-negative-m}}\label{proof-BL-negative-m}

It is well-known that each uncountable compactum contains a Cantor set.
Moreover, 
for any uncountable compactum 
$\mathfrak X\subset S^{N-1} $ 
we can choose a Cantor set 
$\mathfrak X_0\subset \mathfrak X$  
such that $\mathfrak X_0$ is 
cellularly separated in $S^{N-1}$: 
$\mathfrak X_0 = \bigcap\limits_{m=1}^\infty \bigcup\limits_{\ell = 1}^{\ell _m} \overline{V_{m, \ell }}$,
where 
each $V_{m, \ell } \subset S^{N-1}$ is 
homeomorphic to $\mathbb R^{N-1}$,
the closure
$\overline {V_{m, \ell }} $
of each  $V_{m, \ell } $ in $S^{N-1}$
is a topological $(N-1)$-cell,
and 
$\overline { V_{m, \ell  }} \cap \overline { V_{m, n} }= \emptyset $
for each $m\geqslant 1$ and each  $1\leqslant \ell \neq n \leqslant \ell _m$
\cite[I.3, I.4]{Keldysh}.

Take a Krushkal Cantor set $K\subset \mathbb R^N$.
Fix any homeomorphism
$f_0: \mathfrak X_0 \cong K$.
Let  $F : S^{N-1} \to \mathbb R^N$ be any embedding
that extends $f_0$.
(The existence of such extension,
even with the additional property of
being piecewise-linear on $S^{N-1} \setminus \mathfrak X_0$,
is proved using ``horn pulling method''
which is introduced in the works
of L.~Antoine and J.~Alexander;
see \cite[Stat. 4]{Frolkina} for detailed references.
The idea can be found in \cite[Thm. 18.6, 18.7]{Moise}, \cite[Example 2.7.1]{DV}.)
Let $f := F|_{\mathfrak X}$ and $X:= f(\mathfrak X)$.
We have
$$
\dem X = 
\dem f (\mathfrak X) \geqslant
\dem f (\mathfrak X_0) = 
\dem K = N-2  > N-3 = \dem M^{N-3}_N,
$$
which implies 1).
Assertion 2) follows from Corollary \ref{Krushkal-Cantor}.

\subsection{Proof of Theorem \ref{BL-bouquet}}\label{proof-BL-bouquet}

We will construct an embedding
$f: P \rightarrow \mathbb R^N$
which satisfies assumptions of Corollary~\ref{push}.
There exist a point $a\in P$ and a number $\rho >0$
such that $O_{\rho  } (a) \cap P \cong (0, 1)$.
We may moreover assume that
$\overline{O_{\rho  } (a) } \cap P$ is a straight line segment;
denote it $J$ for brevity.

Take 
a 
Blankinship--Ivanov 
Cantor set 
$A\subset \mathbb R^N$ 
(this is a wild Cantor set
constructed for $N\geqslant 4$
using solid tori 
$I^2 \times (S^1)^{N-2}$,
analogously to Antoine's Necklaces 
\cite{Blankinship},
\cite{Ivanov},
\cite{Ivanov-diss},
\cite[Example 4.7.1]{DV}).
We 
may assume 
that $A\subset O_{\rho }(a)$.

As in the proof of Theorem~\ref{BL-negative-m},
construct 
an embedding
$f_0 : J \to O_{\rho } (a)$
such that
$f_0 |_{\partial J} = \id $,
$f_0 (J) \supset A$,
and the restriction of $f_0$ to $J\setminus (f_0)^{-1} (A)$
is piecewise-linear.
Define an embedding 
$f:P\to\mathbb R^N$
by
$$
f(x) = 
\begin{cases}
 x &\text{for } x\in P\setminus J, 
\\
 f_0(x) &\text{for }
x\in J .
\end{cases}
$$
For brevity, denote
$\Sigma := f(J)$ and
$\Pi := \overline {P\setminus J}$.
We have $\Pi \cup \Sigma = f(P) $
and $\Pi \cap \Sigma  = \partial \Sigma $.
Replacing the segment $J$ with a smaller one if necessary, 
we may assume that $\dim \Pi = k$.

The compactum $f(P)$
can not be position-wise embedded into $M^{N-3}_N$
since
$$
\dem f(P)
\geqslant 
\dem A 
= N-2 > N-3 = \dem M^{N-3}_N . 
$$

In the rest of the proof,
we 
show that $f$ satisfies assumptions of Corollary~\ref{push}.

\medskip

\emph{Step 1. Pushing a ``singularity'' $A$ off itself.}

Take a positive number $\alpha $ such that
$ \alpha < d(\Pi , A)$.
Recall that the set $A$ has a special structure: it is an intersection
of sets which are the unions of disjoint solid tori
$I^2\times (S^1)^{N-2}$.
By \cite[p. 171, Remark]{DV} there 
exists an $\alpha $-homeomorphism 
$g_1 : \mathbb R^N \cong\mathbb  R^N$ 
such that
$$
A \cap g_1(A) = \emptyset 
\quad
\text{and}
\quad
\supp g_1 \subset O_{\alpha } (A) .
$$
By the choice of $\alpha $,
we have
$$
\Pi  \cap g_1(A) = \emptyset 
\quad
\text{and}
\quad
 g_1 |_{\Pi } = \id .
$$

\emph{Step 2.  Pushing a ``singularity'' off the arc $\Sigma $.}

Step 2 can be skipped
putting $g_2 := \id $
if $\Sigma \cap g_1(A) = \emptyset $.
If not, do as follows.
Take a positive number $\beta $ such that
$$
 \beta < d(\Pi \cup A , g_1 (A))
 \quad
 \text{ and }
 \quad
 \alpha + \beta < d (\Pi , A).
 $$ 
Let us construct a 
$\beta $-homeomorphism $g_2: \mathbb R^N \cong\mathbb R^N$
such that
$$
\Sigma \cap g_2g_1(A) = \emptyset 
\quad
\text{and}
\quad
\supp g_2 \subset O_{\beta } (\Sigma \cap g_1(A)) .
$$
To prove that the desired $g_2$ exists,
consider 
$L:= \Sigma \setminus O_{\beta }(A)$.
By construction,
$L$ is a finite $1$-dimensional polyhedron.
By the choice of $\beta $,
we have 
$$
\Sigma \cap g_1(A) \cap O_{\beta }(A)
\subset 
g_1(A)\cap O_{\beta }(A) = \emptyset ,
$$
therefore
$$
\Sigma \cap g_1(A) = L\cap g_1(A).
$$
Recall that $\dem g_1(A) = \dem A = N-2$,
and
apply 
the General Position Theorem 
\cite[Thm. 10]{Shtanko-diss}, \cite[Cor. 3.4.7]{DV}
to $L$ and $g_1(A)$
in the $\beta $-neighbourhood of their intersection.
We thus get a desired $g_2$.

Observe that
$$
(\Pi \cup \Sigma )\cap g_2g_1(A)=\emptyset .
$$

\emph{Step 3. Pushing the rest part of $\Sigma $ off the initial compactum $\Pi \cup \Sigma $.}

At Step 3,
we will construct a special self-homeomorphism $g_3$
of $\mathbb R^N$
whose properties include the equality
$$
(\Pi \cup \Sigma ) \cap g_3g_2g_1 (\Sigma ) = \emptyset .
$$
 Take a positive number $\gamma $ such that
$$
 \gamma < d(\Pi \cup \Sigma , g_2g_1 (A))
\quad
\text{ and }
\quad
\alpha +\beta +\gamma < d(\Pi , A).
$$
Take an open neighbourhood $W$ of
$A$ in $\mathbb R^N$
such that
$O_{\gamma } (g_2g_1 (A)) = g_2g_1(W)$.
Denote
$L' := \Sigma \setminus W$.
By construction,
$L'$ is a finite polyhedron with
$\dim L' \leqslant 1$, and
$$
(\Pi \cup\Sigma ) \cap g_2g_1(\Sigma \setminus L')
\subset
(\Pi \cup\Sigma ) \cap g_2g_1(W) =
(\Pi \cup\Sigma ) \cap O_\gamma (g_2g_1(A)) =
\emptyset .
$$
Therefore 
$$
(\Pi \cup \Sigma )\cap g_2g_1(\Sigma ) =
(\Pi \cup \Sigma )\cap g_2g_1(L').
$$
Recall that
$$
\dem (\Pi \cup \Sigma ) = N-2
\quad
\text{ and }
\quad
\dem g_2g_1(L') = \dem L'  = \dim L' \leqslant 1 .
$$
By the General Position Theorem
\cite[Thm. 10]{Shtanko-diss}, \cite[Cor. 3.4.7]{DV},
there exists a
$\gamma $-homeomorphism $g_3: \mathbb R^N \cong\mathbb R^N$
such that
\begin{equation}
\label{3-1}
(\Pi \cup \Sigma )\cap g_3g_2g_1(L') = \emptyset 
\quad
\text{and}
\quad
\supp g_3 \subset O_{\gamma } ((\Pi \cup \Sigma )\cap g_2g_1(L')) 
\end{equation}
By construction, 
$g_3$ restricted on $ O_{\gamma } (g_2g_1(A)) = g_2g_1(W)$
is the identity map,
hence
\begin{equation}
\label{3-2}
(\Pi \cup \Sigma )\cap g_3g_2g_1(\Sigma \cap W)=
(\Pi \cup \Sigma )\cap g_2g_1( W) = \emptyset  
\end{equation}
From $(\ref{3-1})$ and $(\ref{3-2})$
we get the desired property:
$$
(\Pi \cup \Sigma )\cap g_3g_2g_1(\Sigma  ) = \emptyset .
$$
In addition, the inequality $\alpha +\beta +\gamma < d(\Pi , A)$
implies
$$
A\cap g_3g_2g_1(\Pi ) = \emptyset  .
$$

\emph{Step 4.  Pushing the rest part of the compactum $\Pi \cup \Sigma $ off itself.}

Take a positive number $\delta $ such that
$$
 \delta < \min \{ 
d(\Pi \cup \Sigma , g_3g_2g_1 (\Sigma ));
\quad
d(A , g_3g_2g_1 (\Pi  )) 
\}
.
$$
Denote
$L'' := \Sigma \setminus O_{\delta }(A) $.
Again, 
$L''$ is a finite polyhedron with
$\dim L'' \leqslant 1$.
Observe that
$$
\dem (\Pi \cup L'' ) = 
k
\quad
\text{ and }
\quad
\dem g_3g_2g_1(\Pi ) = \dem \Pi = k.
$$
By assumption, $2k+1\leqslant N$.
Applying the General Position Theorem
to 
$\Pi \cup L''$
and 
$g_3g_2g_1(\Pi ) $,
we find
a
$\delta $-homeomorphism $g_4: \mathbb R^N \cong\mathbb R^N$
such that
$$
(\Pi \cup L'' )\cap g_4g_3g_2g_1(\Pi ) = \emptyset 
\quad
\text{and}
\quad
\supp g_4 \subset O_{\delta } 
((\Pi \cup L'' )\cap g_3g_2g_1(\Pi )) .
$$
By construction,
$O_{\delta }(A) \cap g_4g_3g_2g_1(\Pi )= \emptyset $.
Hence
\begin{equation}
\label{Pi}
(\Pi \cup \Sigma )\cap g_4g_3g_2g_1(\Pi ) = \emptyset  
\end{equation}
Observe that
$g_4$ restricted on
$O_{\delta } (g_3g_2g_1 (\Sigma ))$
is the identity map, consequently
\begin{equation}
\label{Sigma}
(\Pi \cup \Sigma )\cap g_4g_3g_2g_1(\Sigma ) = 
(\Pi \cup \Sigma )\cap g_3g_2g_1(\Sigma ) = \emptyset  
\end{equation}
From 
(\ref{Pi})
and
(\ref{Sigma})
we finally obtain
$$
(\Pi \cup \Sigma )\cap g_4g_3g_2g_1(\Pi \cup \Sigma ) 
= \emptyset .
$$

The composition
$h:= g_4g_3g_2g_1$
is the desired homeomorphism.
Indeed: 
$\supp h \subset O_{\alpha +\beta +\gamma +\delta } (f(P))$,
and $h$ moves each point 
of $\mathbb R^N$ no further than
$\alpha +\beta +\gamma +\delta $
(which can be made smaller than an arbitrary given
positive number
$\varepsilon $).

\subsection{Proof of Theorem \ref{BL-disconnected}}\label{proof-BL-disconnected}

By assumption, $N-k\geqslant 3$.
Take an embedding 
$f: \mathcal C\times\mathcal C \rightarrow \mathbb R^{N-k}$
such that all
$A_t  :=  f(\mathcal C \times \{ t\})$, where $t\in\mathcal C$,
are 
wild Cantor sets 
embedded in $\mathbb R^{N-k}$
equivalenly to each other.
For example, we may take Antoine--Blankinship--Ivanov sets
\cite[Remark on p. 171 and Exercise 4.8.1]{DV}.
Applying the idea from \cite[p. 479]{Bryant},
define a new embedding $F$ by multiplying with $I^k$:
$$
I^k\times (\mathcal C\times\mathcal C )
\stackrel{\id\times f}{\rightarrow }
I^k \times \mathbb R^{N-k}
\subset 
\mathbb R^k \times \mathbb R^{N-k}
=
\mathbb R^{N} .
$$
For each $t\in\mathcal C$,
denote
$X_t := 
F( I^k\times \mathcal C\times\{ t\} ) = I^k\times A _t \subset \mathbb R^N $.
All $X_t$'s are embedded into $\mathbb R^{N}$
equivalently to each other.

Recall that a wild Cantor set in $\mathbb R^N$
can not be locally $1$-co-connected;
see \cite{Homma} or \cite[Thm. 5.1]{Bing-tame} for $N=3$,
\cite{McMillan-taming} for $N\geqslant 5$,
and \cite{Shtanko-4dim} or \cite[p.~5]{BDVW} for $N=4$.
(A closed subset $X$ of $\mathbb R^d$
is called \emph{locally
$1$-co-connected}, or briefly $1$-LCC
if
for each $x\in \mathbb X$,
any neighbourhood $U$ of $x$ in $\mathbb R^d$
contains a smaller neighbourhood $V$ of $x$
such that
every map
$\gamma : S^1 = \partial (I^2) \to V\setminus X$
can be extended to a map
$\Gamma : I^2 \to U\setminus X$.
Refer to \cite{Stanko1970}, \cite{Edwards}, \cite[3.4]{DV}, \cite{Cernavsky}
for details.)

Hence for any $t\in\mathcal C$,
the set $A_t$
is not 
$1$-LCC
in 
$\mathbb R^{N-k} $.                                                                                                                                                                                                                                       
Consequently,
$X_t$ is not $1$-LCC in
$\mathbb R^{N}$.
By
\cite[proof of Thm.~2]{Bryant1969} or
\cite{Stanko1970}
(one may also refer to \cite[Thm. 1.4]{Edwards}
or \cite[Thm. 3.4.11]{DV}),
$X_t$ can not be position-wise 
embedded into $M^{N-3}_{N}$.

\section{Appendix}\label{App}

The construction of Baker and Laidacker
includes careful consideration of metric properties, providing $\mathfrak c$
congruent compacta.
In this section, we give a short topological proof
for a weaker statement (Corollary \ref{BL-top}).
We also state several consequences,
including a generalization of
results obtained by R.B.~Sher and E.R.~Apodaca.

\begin{cor}\label{BL-top}
For each $k\geqslant 0$,
there exists an embedding
$\psi : M_{2k+1}^k \times\mathcal C \rightarrow \mathbb R^{2k+1}$
such that for each $t\in\mathcal C$ the image
$\psi (M_{2k+1}^k\times \{ t\})$ 
is ambiently homeomorphic to the standard Menger compactum
$M_{2k+1}^k \subset \mathbb R^{2k+1}$.
\end{cor}

\begin{pf}
The standard Menger compactum 
$M_{2k+1}^k $
can be taken off 
itself by a small self-homeomorphism of
$\mathbb R^{2k+1}$
with support
in an arbitrarily small neighbourhood of $M^k_{2k+1}$
\cite[Prop. 1]{Stanko1971},
\cite[Prop. 8, Thm. 10]{Shtanko-diss},
\cite[Cor. 3.4.7, Thm. 3.5.1]{DV}.
The conditions of Corollary~\ref{push} are thus satisfied.
\end{pf}

The Baker--Laidacker theorem
implies the following.
Let $\mathcal H$ be a 
family 
of 
compacta whose dimensions do not exceed $k$,
and
$| \mathcal H | \leqslant\mathfrak c$;
then 
it is possible to embed all members of $\mathcal H$
in 
$\mathbb R^{2k+1}$
simultaneously,
so that the images are pairwise disjoint
\cite[Cor.~1]{BL}.
Observing
that there are exactly $\mathfrak c$
compact subsets of $\mathbb R^{2k+1}$,
we immediately get:

\begin{cor}\label{BL-simult-n}
For each $k\geqslant 0$,
there exists a family $\mathcal H$
of pairwise disjoint compacta 
in $\mathbb R^{2k+1}$
such that 
for each compactum $X$ with $\dim X \leqslant k$
the family $\mathcal H$
contains $\mathfrak c$
elements homeomorphic to $X$.
\end{cor}

Assuming the
Continuum Hypothesis,
R.B.~Sher 
proved
\cite[Thm. 3]{Sher1969}:
there is a collection
$\mathcal H$ of pairwise disjoint
arcs in
$\mathbb R^3$ 
such that if $A\subset \mathbb R^3$ is an arc 
whose wild set is a compact $0$-dimensional set,
then there is an $A' \in \mathcal H$
which is embedded in $\mathbb R^3$
equivalently to $A$.

Actually, 
without the Continuum Hypothesis we show that
$\mathbb R^3$ is ``a universal storage'' of all possible knots and arcs
(this also implies \cite[Thm. 4,~5]{Apodaca}):

\begin{cor}\label{knots}
There exists a family $\mathcal H$
of 
pairwise disjoint
compacta in~$\mathbb R^3$ 
with the property:
for any one-dimensional compactum
 $X\subset \mathbb R^3$ 
that can be represented as the union
 of 
no more than countably many ANR-sets,
we have 
$$| \{ Y \in \mathcal H \ | \ 
(\mathbb R^3 , Y) \cong (\mathbb R^3 , X) \} | = 
\mathfrak c .
$$
\end{cor}

Recall that a metric space $X$ is
an ANR-space (an absolute neighbourhood retract)
if 
for any metric space $Y$, any closed subset
$Y_0\subset Y$, and any continuous map
$f: Y_0 \to X$
there exists a continuous extension $F:U(Y_0)\to X$
of $f$ to an open neighbourhood of $Y_0$.
By the Borsuk theorem, any polyhedron is an ANR.

\begin{pf-cor-knots}
Any $X$ which 
satisfies our assumptions 
can be position-wise embedded into $M^1_3$
\cite[Satz 4]{Bothe1964univ}.
The set of all compact subsets of $\mathbb R^3$
has cardinality
$\mathfrak c$, 
hence the desired statement follows from \cite[Thm.~1]{BL}
or from Corollary \ref{BL-top}.
\end{pf-cor-knots}

Similarly to Corollary \ref{knots}, we get

\begin{cor}\label{Cantor-sets}
There exists a family $\mathcal H$
of 
pairwise disjoint
compacta in~$\mathbb R^3$ 
such that
for any zero-dimensional compactum
 $X\subset \mathbb R^3$ 
there is an $Y \in \mathcal H$
which is embedded in $\mathbb R^3$
equivalently to $X$.
\end{cor}

\begin{pf}
Each zero-dimensional compactum
$X\subset \mathbb R^3$ 
can be position-wise embedded into $M^1_3$;
this can be derived from
\cite[Thm.~1]{Sher1969}
using the idea from \cite[Proof of Satz~4]{Bothe1964univ}.
Alternatively, 
we may directly apply 
\cite[Satz~4]{Bothe1964univ} 
if we
recall the Denjoy--Riesz theorem:
for each $N\geqslant 1$
and each zero-dimensional
compactum
$X\subset \mathbb R^N$
there exists an arc $\widehat X\subset \mathbb R^N$
with
$X\subset \widehat X$.
\end{pf}

\section*{Final remarks}

\begin{rmk}
It remains an open question whether
there exists a $1$-dimensional compactum $X\subset \mathbb R^3$ such that:
$X$ can not be position-wise embedded into
$M^1_3$, and
$\mathbb R^3$ contains 
an uncountable family of pairwise disjoint compacta
embedded equivalently to $X$.
\end{rmk}

\begin{rmk}
For each
$N\geqslant 4$,
D.G.~Wright 
described
cellular arcs $L_1, L_2$ in $\mathbb R^N$
such that $L_1$ cannot be slipped off $L_2$,
and consequently
$X=L_1\cup L_2$
is sticky
\cite{Wright-sticky}.
This could possibly be used to construct new
``negative'' examples for the question of
B.J.~Baker and M.~Laidacker.
\end{rmk}

\end{document}